\documentclass[12pt]{article}
\usepackage{hyperref}
\usepackage{amssymb, amsmath}
\usepackage{graphicx}
\usepackage[mathscr]{eucal}

\newtheorem{thrm}{Theorem}

\newtheorem{lemma}{Lemma}

\begin{document}
\title{The spectrum of the Laplacian in a domain bounded by 
a flexible polyhedron in $\mathbb R^d$ does not always 
remain unaltered during the flex}
\author{Victor Alexandrov}
\date{March 16, 2019}
\maketitle

\begin{abstract}
Being motivated by the theory of flexible polyhedra,
we study the Dirichlet and Neumann eigenvalues 
for the Laplace operator in special bounded domains
of Euclidean $d$-space.
The boundary of such a domain is an embedded simplicial
complex which allows a continuous deformation (a flex),
under which each simplex of the complex moves as a solid
body and the change in the spatial shape of the domain is achieved
through a change of the dihedral angles only.
The main result of this article is that
both the Dirichlet and Neumann spectra of the Laplace operator
in such a domain do not necessarily remain unaltered during the flex
of its boundary.
\end{abstract}

\section{Introduction}\label{s1}

In this article, a \textit{polyhedron}
is either a continuous map $f: K\to\mathbb{R}^d$ of a connected 
$(d-1)$-dimensional simplicial complex $K$ which is affine on every simplex 
or the image  $f(K)\subset\mathbb{R}^d$ of $K$ under the action of $f$.
A polyhedron is called \textit{embedded} if $f$ is injective.
A polyhedron $P_0=f(K)$ is called \textit{flexible} if its spatial shape can be
changed continuously only by changing its dihedral angles attached to its
$(d-2)$-dimensional faces, i.\,e. by means of a continuous deformation,
which does not change the dimensions of the faces.
In other words, $P_0$ is said to be flexible if
there exists a continuous family $\{P_t \}_{t\in[0,1]}$ of
polyhedra such that, for every $t\in(0,1]$, $P_0$ and $P_t$ are 
combinatorially equivalent to each other and every two corresponding faces of
$P_0$ and $P_t$ are congruent, while $P_0$ and $P_t$ themselves are not congruent.
Wherein, the family $\{P_t \} _{t\in[0,1]}$ is called a \textit{flex} of $P_0$.

The above definitions are standard in the theory of flexible polyhedra.
However, in this article, we will also use the following notation, which is 
not commonly used: by $[\![P_0]\!]$ we denote the bounded domain in $\mathbb{R}^d$,
whose boundary is an embedded boundary-free polyhedron $P_0\subset\mathbb{R}^d$.

The fact that embedded boundary-free flexible polyhedra do exist in $\mathbb{R}^3$  
is non-trivial and was established only in 1977 by Robert Connelly \cite{Co76}
(see also \cite{Co80}, \cite {Ku79}).
For $d\geqslant 4$, the question about the existence of an embedded boundary-free 
flexible polyhedron in $\mathbb{R}^d$ remains open.

Over the last 40 years, through the efforts of many geometers,
it was shown that flexible polyhedra (not necessarily embedded) have
a number of remarkable properties. Below we mention some of them:

(i) \textit{In $\mathbb{R}^3$, there exists an embedded sphere-homeomorphic flexible
polyhedron with nine vertices only.}
This polyhedron was built by Klaus Steffen around 1980 and bears his name.
Although Steffen published no text about his polyhedron,
its description can be found in many books and
articles (see, e.\,g., \cite[Section 23.2.3]{DO07} or \cite{Al10}).
The question about the existence of an embedded boundary-free flexible polyhedron in
$\mathbb{R}^3$ with eight vertices remains open, see \cite{Ma08}.

(ii) \textit{In $\mathbb{R}^3$, there exist embedded boundary-free flexible
polyhedra of an arbitrary genus, both orientable and non-orientable.} 
Explicit examples of such polyhedra can be found in \cite{Sh15}.

(iii) With every oriented boundary-free polyhedron in $\mathbb{R}^d$
one can associate a quantity called its \textit{oriented volume}.
For an embedded boundary-free polyhedron $P_0$ in $\mathbb{R}^d$, its oriented
volume coincides, up to a sign, with the $d$-dimensional volume of the domain $[\![P_0]\!]$.
So, \textit{for any $d\geqslant 3$ and any oriented boundary-free flexible polyhedron
$P_0$ in $\mathbb {R}^d$, the oriented volume of $P_0$ remains constant during the flex.}
For several decades, this statement was known as the Bellows Conjecture.
In $\mathbb{R}^3$, it was first proved by I.Kh. Sabitov in 1995--1996 in 
\cite{Sa95}, \cite {Sa96} (see also \cite {Sa98}).
Another proof, which is also limited to the case $d=3$, was published in 1997 in \cite {CSW97}
(see also an expository paper \cite{Sc04}).
For $d\geqslant 4$, the Bellows Conjecture was proved by A.A. Gaifullin in 2014 
in \cite{Ga14a} and \cite{Ga14b} (see also an expository paper \cite{Ga16}).

(iv) With every oriented polyhedron in $\mathbb{R}^d$
one can associate a quantity called its \textit{integral mean curvature}, 
see, e.\,g., \cite[Chapter 13]{Sa76}.
In the case of an oriented polyhedron in $\mathbb{R}^3$, its
integral mean curvature is given by
$\frac12\sum_{i}(\pi-\varphi_i)\ell_i$,
where the sum is taken over all edges of the polyhedron,
$\ell_i$ is the length of the $i$-th edge and $\varphi_i$ is
the dihedral angle along the edge.
So, \textit{for every $d\geqslant 3$, the integral mean
curvature of an oriented flexible polyhedron
$P_0\subset\mathbb{R}^d$ remains constant during the flex.}
This statement was proved by Ralph Alexander in 1985
in \cite{Al85}.

(v) \textit{For every $d\geqslant 3$, if an embedded
boundary-free polyhedron $P_1\subset\mathbb{R}^d$ is obtained from an embedded
boundary-free polyhedron $P_2\subset\mathbb{R}^d$ by a flex
then $[\![P_1]\!]$ and $[\![P_2]\!]$ have the same Dehn invariants.}
This statement was proved by A.A. Gaifullin and L.S. Ignashchenko in 2017 in \cite{GI17}.
Since it is known that domains with polyhedral boundaries in $\mathbb{R}^3$
are scissors congruent if and only if they have the same volume and
the same Dehn invariants, then from the property under discussion
it follows immediately that \textit{if an embedded boundary-free
polyhedron $P_1\subset\mathbb{R}^3$ is obtained from an embedded
boundary-free polyhedron $P_2\subset\mathbb{R}^3$ by a flex
then the domains $[\![P_1]\!]\subset\mathbb{R}^3$ and
$[\![P_2]\!]\subset\mathbb{R}^3$ are scissors congruent.}
Since 1980s, the latter statement was known as the Strong Bellows Conjecture.

Note that, for every $d\geqslant 3$, the notion of the flexible polyhedron
can be defined in any $d$-dimensional space of constant curvature (i.\,e., not 
only in Euclidean space, but also in spherical and hyperbolic spaces),
as well as in Minkowski space.
The reader, interested in properties of flexible polyhedra in these spaces, which are
similar to the properties (i)--(v), is referred to \cite{Al03}, \cite{Ga14c}--\cite{Ga15b}, 
\cite{Ga17}, \cite{St00}, and \cite{St06}.

Being motivated by the properties (i)--(v), we would like to find new invariants of 
flexible polyhedra in $\mathbb{R}^d$, $d\geqslant 3$, which are preserved during the flex.
In our opinion, it is natural to try the Dirichlet and Neumann eigenvalues of the  
Laplace equation in the domain $[\![P_0]\!]$ on the role of such invariants.
The statement that the Dirichlet and Neumann spectra of the Laplacian are not altered 
during the flex seems natural  to us because  it agrees with the Weyl law on the asymptotics of 
eigenvalues of the Laplacian. 
In fact, the coefficients of the first and second terms of the Weyl asymptotics are expressed 
in terms of the volume and the surface area of the boundary of the domain 
and, thus, remain unaltered during the flex (see the formula (\ref {eq3}) in Section~\ref{s2},
where we recall the Weyl law in more detail).

The main result of this article is the following theorem showing that the above conjecture on 
the invariance of the Dirichlet and Neumann spectra of the Laplacian 
during the flex of the boundary of a domain is false:

\begin{thrm}\label{thrm1}
For every $d\geqslant 3$, $\varepsilon >0$, and every embedded
flexible polyhedron $P_0\subset\mathbb{R}^d$ there is an
embedded flexible polyhedron $\widetilde{P}_0\subset\mathbb{R}^d$ 
and its flex $\mathscr{\widetilde{F}}=\{\widetilde{P}_s\}_{s\in [0,1)}$ such that

{\rm ($\alpha_d$)}
the combinatorial structure of $\widetilde {P}_0$ is a subdivision of the combinatorial 
structure of $P_0$;

{\rm ($\beta_d$)}
the Hausdorff distance between the sets $\widetilde{P}_0$ and $P_0$ is less than~$\varepsilon$;

{\rm ($\gamma_d$)}
both Dirichlet and Neumann spectra of the $d$-dimensional Laplacian in the domain 
$[\![\widetilde{P}_s]\!]\subset\mathbb{R}^d$ do not remain unaltered when $s$ changes in the interval $[0,1)$.
\end{thrm}

\textbf{Remark 1} \ 
For $d=2$, the following statement, similar to Theorem \ref{thrm1}, holds true:
\textit{For every $\varepsilon>0$ and every closed embedded polygon $P_0\subset\mathbb{R}^2$, 
there is a closed embedded polygon $\widetilde{P}_0\subset\mathbb{R}^2$ and its flex 
$\mathscr{\widetilde{F}}=\{\widetilde{P}_s\}_{s\in [0,1)}$ such that}
\par
($\alpha_2$)
\textit{the number of vertices of $\widetilde{P}_0$
exceeds the number of vertices of $P_0$ by no more than four};
\par
($\beta_2$)
\textit{the Hausdorff distance between the polygons $\widetilde{P}_0$
and $P_0$ is less than~$\varepsilon$};
\par
($\gamma_2$)
\textit{
both Dirichlet and Neumann spectra of the $2$-dimensional Laplacian in the domain 
$[\![\widetilde{P}_s]\!]\subset\mathbb{R}^2$ 
do not remain unaltered when $s$ changes in the interval $[0,1)$.}

This statement does not belong to the theory of flexible polyhedra.
For this reason, we prefer to formulate it separately from Theorem \ref{thrm1}.
The proof of the statement under discussion is left to the reader.
It can be obtained by simplifying the proof of Theorem \ref{thrm1} which is given in Section \ref{s4}.

\textbf{Remark 2} \ 
It would be interesting to obtain an analogue of Theorem \ref{thrm1}
for spherical and hyperbolic spaces of dimension $d\geqslant 2$.
We cannot do this in this paper, because our proof of Theorem \ref{thrm1}, presented in 
Section~\ref {s4}, relies on an asymptotic formula for the eigenvalues of the Laplace operator
given below in Lemma \ref {l1}, while we are not aware about any analogue of Lemma \ref{l1} for spherical or
hyperbolic spaces.

\section{Weyl's law and Fedosov's asymptotic formula}\label{s2}

Let $d\geqslant 2$, $\Omega$ be a bounded domain in $\mathbb{R}^d$, and
$\Delta=\sum_{i=1}^d {\partial^2}/{\partial x_i^2}$ be the Laplace operator in $\mathbb{R}^d$.
Consider the following boundary problems:
\begin{equation}\label{eq1}
\left\{
\begin{array}{rll}
\Delta u&=-\nu^2 u &\quad \mbox{in}\ \Omega,\\
u\vert_{\partial\Omega}&=0& 
\end{array}
\right. 
\end{equation}
and
\begin{equation}\label{eq2}
\left\{
\begin{array}{ll}
\Delta u&=-\nu^2 u \quad \mbox{in}\ \Omega,\\
\frac{\partial u}{\partial\boldsymbol{n}}\vert_{\partial\Omega}&=0,
\end{array}
\right. 
\end{equation}
where ${\partial u}/{\partial\boldsymbol{n}}$ is the directional 
derivative of the function $u$ in the direction $\boldsymbol{n}$ 
and $\boldsymbol{n}$ is the exterior unit normal to $\partial\Omega$. 

As usual, we call the number $\nu^2$ for which there exists a nonzero
the solution $u$ of the problem (\ref{eq1}) (respectively, (\ref{eq2})) 
the \textit{Dirichlet eigenvalue}  (respectively, the \textit{Neumann eigenvalue})
of the Laplace operator in $\Omega$,
For each of the problems (\ref{eq1}) and (\ref{eq2}), we denote by
$\mathscr{N}(k)$ the number of eigenvalues, which do not exceed $k^2$ 
(repeating each eigenvalue according to its multiplicity).
We call $\mathscr{N}(k)$ the \textit{eigenvalue counting function}
of the corresponding problem.

Under certain assumptions about the boundary $\partial\Omega$ of a
bounded domain $\Omega\subset\mathbb{R}^d$,
the following asymptotic formula holds true for $k\to\infty$:
\begin{equation}\label{eq3}
\mathscr{N}(k)=\frac{\mbox{vol}_d(\Omega)}{\Gamma\bigl(\frac{d+2}{2}\bigr)}
\biggl(\frac{k}{2\sqrt{\pi}}\biggr)^d\mp
     \frac{\mbox{vol}_{d-1}(\partial\Omega)}{4\Gamma\bigl(\frac{d+1}{2}\bigr)}
\biggl(\frac{k}{2\sqrt{\pi}}\biggr)^{d-1}+o(k^{d-1}).
\end{equation}
Here and below, $\mbox{vol}_p$ denotes the $p$-dimensional volume of a set and
$\Gamma$ denotes the Euler gamma function.
The minus sign corresponds to the problem (\ref{eq1}),
while the plus sign corresponds to the problem (\ref{eq2}).

The formula (\ref{eq3}) is known as Weyl's law or the Weyl asymptotic formula.
The reader interested in more details about this formula and its influence on mathematics and physics
is referred to \cite{ANPS}, \cite {Iv16}, and literature mentioned there.

If the boundary $\partial\Omega$ of the domain $\Omega\subset\mathbb{R}^d$
is a flexible polyhedron, then the coefficients of $k^d$ and
$k^{d-1}$ in the formula (\ref{eq3}) remain constant during the flex of $\partial\Omega$.
This observation is our main argument supporting the conjecture on 
the invariance of the Dirichlet and Neumann spectra of the Laplacian 
formulated in Section~\ref{s1}.

Direct calculation of the eigenvalues of the Laplacian for domains bounded by flexible polyhedra
is hardly possible due to the complex geometry of such domains.
The use of the formula (\ref{eq3}) is also hardly possible because finding every new term in
this asymptotics is a difficult problem (see, e.\,g., \cite{BS88}, \cite{Sm81}, \cite{Wa05}) 
and the nonsmoothness of the boundary of the domain gives rise to 
additional difficulties (see, e.\,g., \cite{GL17}, \cite{NS05}).
Nevertheless, we refute the conjecture under discussion.
For this we use the following version of Weyl's law for the Riesz means of the eigenvalues of 
the Laplacian in a domain bounded by a polyhedron in Euclidean $d$-space, $d\geqslant 2$:

\begin{lemma}[B.V. Fedosov \cite{Fe64}]\label{l1}
Let $d\geqslant 2$, $0\leqslant p\leqslant d-1$, and let a bounded domain $D\subset\mathbb{R}^d$ 
be such that its boundary $\partial D$ is a polyhedron. 
Let $\{ F^{d-2}_i\}_i$ be the set of all 
$(d-2)$-dimensional faces of $\partial D$, and 
let $\varphi_i$ stand for the value of the dihedral angle of $D$ at $F^{d-2}_i$.  
Then the following asymptotic formula, involving the eigenvalue counting function
$\mathscr{N}(k)$, holds true for each of the problems $(\ref{eq1})$ and $(\ref{eq2}):$
\begin{equation*}
\frac{1}{\Gamma (p+1)} \int\limits_0^k(k-t)^p\,d\mathscr{N}(t)=
\sum\limits_{l=1}^d a_l\frac{\Gamma(l+1)}{\Gamma(p+l+1)}k^{p+l}+
O(k^{d-1})
\end{equation*}
as $k\to\infty$. Here
\begin{eqnarray}
a_d&=&\frac{\mbox{\rm{vol}}_d(D)}{2^d \pi^{d/2} 
\Gamma\bigl(\frac{d}{2}+1\bigr)},\nonumber\\
a_{d-1}&=&\mp
\frac{\mbox{\rm{vol}}_{d-1}(\partial D)}{2^{d+1}\pi^{(d-1)/2}
\Gamma\bigl(\frac{d+1}{2}\bigr)},\label{eq4}\\
a_{d-2}&=&\frac{1}{2^{d+1}\pi^{d/2}\Gamma\bigl(\frac{d}{2}\bigr)}
\sum_i\frac{\varphi^2_i-\pi^2}{3\varphi_i}
\mbox{\rm{vol}}_{d-2}\bigl(F^{d-2}_i\bigr).\label{eq5}
\end{eqnarray}
In the formula $(\ref{eq4})$, the minus sign corresponds to the problem $(\ref{eq1})$,
while the plus sign corresponds to the problem $(\ref{eq2})$.
\end{lemma}

Lemma \ref{l1} was proved by B.V. Fedosov in 1964 in \cite{Fe64}, where he used the same method 
as in \cite{Fe63}, where he proved Lemma \ref{l1} for the case $d=2$.
From our point of view, the article \cite{Fe63} is an inalienable part of the article \cite{Fe64}, 
because the latter does not contain details which are common for the cases $d=2$ and $d\geqslant 3$.

Note also that, for the case $d=2$, the formula (\ref{eq5}) is discussed in \cite{MR15}.

In Section \ref{s4} we show that Theorem \ref{thrm1} can be deduced from the formula (\ref{eq5}).

\section{Special flexible polyhedron}\label{s3}

To prove Theorem \ref{thrm1}, we need not only the formula (\ref{eq5}), but also some
special flexible polyhedron $\widetilde{P}_0 $, whose existence is established by the following lemma:

\begin{lemma}\label{l2}
Let $d\geqslant 3$, $\varepsilon>0$, and let $P_0$ be an embedded flexible polyhedron in $\mathbb{R}^d$.
Then there is an embedded flexible polyhedron $\widetilde{P}_0\subset\mathbb{R}^d$ 
and its flex $\mathscr{\widetilde{F}}=\{\widetilde{P}_s\}_{s\in [0,1)}$ such that

{\rm ($\widetilde{\alpha}$)} 
the combinatorial structure of $\widetilde {P}_0$ is a subdivision of the combinatorial 
structure of $P_0$;

{\rm ($\widetilde{\beta}$)}
the Hausdorff distance between $\widetilde{P}_0$ and $P_0$ is less than~$\varepsilon$;

{\rm ($\widetilde{\gamma}$)}
for every $s\in(0,1)$, $\widetilde{P}_s$ is an embedded polyhedron; moreover, there is a
$(d-2)$-dimensional face $\widetilde{x}_0$ of $\widetilde{P}_0$ such that the $(d-2)$-dimensional
face $\widetilde{x}_s$ of $\widetilde{P}_s$, which corresponds to $\widetilde{x}_0$ according to the property 
$(\widetilde{\alpha})$, possesses the following properties:

{\quad}\quad{\rm ($\widetilde{\gamma}_1$)} 
the value of the interior dihedral angle $\varphi(\widetilde{x}_s)$ of $\widetilde{P}_s$
at $\widetilde{x}_s$ tends to zero as $s\to 1$; 

{\quad}\quad{\rm ($\widetilde{\gamma}_2$)} 
there exist constants $\widetilde{m}$ and $\widetilde{M}$ such that
$0<\widetilde{m}<\widetilde{M}<2\pi$ and, for every $(d-2)$-dimensional face
$\widetilde{y}_0$ of $\widetilde{P}_0$, which is different from
$\widetilde{x}_0$, and for every $s\in [0,1)$ the inequality
$\widetilde{m}<\varphi(\widetilde{y}_s)<\widetilde{M}$
holds true.  
\end{lemma}

\textbf{Proof} \ Since the polyhedron $P_0$ is assumed to be flexible, there exists a continuous 
family $\mathscr{F}=\{P_t\}_{t\in [0,1]}$ of polyhedra such that, for every $t\in (0,1]$, 
$P_0$ and $P_t$ are combinatorially equivalent, every two of their corresponding faces are congruent
to each other, but the polytopes $P_0$ and $P_t$ themselves are not congruent.

In the process of proving Lemma \ref{l2}, we will construct the polyhedron $\widetilde{P}_0$ and 
the family $\mathscr{\widetilde{F}}=\{\widetilde{P}_s\}_{s\in [0,1)}$ by modifying $P_0$ and 
$\mathscr{F}$ in a special way. We will implement this modification in several steps.

\textit{Step 1.} 
Since $P_0$ is embedded, none of it dihedral angles is equal to 0 or $2\pi$.
Hence, there exist constants $m$ and $M$ such that $0<m<M<2\pi$ and,
for every $(d-2)$-dimensional face $x_0^{(d-2)}$ of $P_0$,
the dihedral angle $\varphi(x_0^{(d-2)})$ of $P_0$ at $x_0^{(d-2)}$ 
satisfies the following inequalities: $m<\varphi(x_0^{(d-2)})<M$.

Since the family $\mathscr{F}=\{P_t\}_{t\in [0,1]}$ is continuous, 
there exists $\omega>0$ such that, for all $t\in [0,\omega]$,
the following properties are valid:

$(\omega_1)$ $P_t$ is an embedded polyhedron; and 

$(\omega_2)$ for every $(d-2)$-dimensional face $x_t^{(d-2)}$ of $P_t$,
the following inequality holds true: $0<m/2<\varphi(x_t^{(d-2)})<M/2+\pi<2\pi$.

Substituting $t\mapsto t/\omega$, we may assume without loss of generality that 
$\omega=1$, i.\,e., that the properties $(\omega_1)$ and $(\omega_2)$ are valid for all $t\in [0,1]$.

\textit{Step 2.} 
Let's choose a $(d-2)$-dimensional face $X_0^{(d-2)}$ of $P_0$ such that, for some $\delta>0$, 
the dihedral angle $\varphi(X_t^{(d-2)})$ is not constant on the interval $[0,\delta]$.
Such a face $X_0^{(d-2)}$ does exist because $P_0$ it is assumed to be flexible.
By definition, put $T=\{t\in(0,\delta) \vert \varphi(X_t^{(d-2)})>\varphi(X_0^{(d-2)})\}$.
Since the function $t\mapsto \varphi(X_t^{(d-2)})$ is continuous, $T$ is an open subset of $\mathbb{R}$. 

By definition, put $t_*=0$, if $t=0$ is a limit point of $T$. If $t=0$ is not a limit point of $T$, 
choose $t_*>0$ such that $t_*$ does not lie in $T$ and is so small that the Hausdorff distance between $P_0$ 
and $P_{t_*}$ is less than  $\varepsilon/3$. In the both cases, put by definition
\begin{equation}\label{eq6}
\varphi_*=\varphi(X_{t_*}^{(d-2)})
\end{equation}
and $T_*=\{t\in(0,\delta) \vert \varphi(X_t^{(d-2)})>\varphi_*\}$.

Choose an interval $(t_1, t_2)\subset T_*$ such that only one of the points $t_1$, $t_2$ belongs to 
the boundary of $T_*$, i.\,e., such that the equality $\varphi(X_t^{(d-2)})=\varphi_*$ is fulfilled 
for $t=t_1$ and is violated for $t=t_2$, or it is fulfilled for $t=t_2$ and is violated for $t=t_1$.

Let us introduce a new parameter $s$ (which is linearly expressed in terms of the parameter $t$)
such that the segment  $[t_1,t_2]$ (within which changes $t$) corresponds to the segment $[0,1]$
(within which changes $s$) and the equality $\varphi(X_s^{(d-2)})=\varphi_*$ holds true for $s=1$.
For the corresponding values of $t$ and $s$, we put $P'_s=P_t$.
Denote by $\mathscr{F}'=\{P'_s\}_{s\in[0,1]}$ the family of polyhedra obtained in this way.

\textit{Step 3.} 
Denote by ${X'_0}^{(d-2)}$ the $(d-2)$-dimensional face of $P'_0$, which corresponds to 
the $(d-2)$-dimensional face ${X_0}^{(d-2)}$ of $P_0$ which was chosen in Step 2.
Then the dihedral angle $\varphi({X'_s}^{(d-2)})$ is not constant for $s\in [0,1]$, 
moreover $\varphi({X'_s}^{(d-2)})>\varphi_*$ for all $s\in[0,1)$
and $\varphi({X'_s}^{(d-2)})=\varphi_*$ for $s=1$.

Denote by ${Y'_0}^{(d-1)}$ and ${Z'_0}^{(d-1)}$ the $(d-1)$-dimensional faces of $P'_0$
incident to  ${X'_0}^{(d-2)}$. 
Let $\xi'$ be a triangulation of ${X'_0}^{(d-2)}$, $\eta'$ be a triangulation of ${Y'_0}^{(d-1)}$ and
$\zeta'$ be a triangulation of ${Z'_0}^{(d-1)}$ such that the following properties hold true:

$(\xi_1)$ the restriction of $\eta'$ on ${X'_0}^{(d-2)}$ coincides with $\xi'$;

$(\xi_2)$ the restriction of $\zeta'$ on ${X'_0}^{(d-2)}$ coincides with $\xi'$;

$(\xi_3)$ there exist simplicies ${y'_0}^{(d-1)}\in \eta'$ and
${z'_0}^{(d-1)}\in \zeta'$ such that
${y'_0}^{(d-1)}\cap {X'_0}^{(d-2)}={z'_0}^{(d-1)}\cap {X'_0}^{(d-2)}$
and the Hausdorff distance between the sets 
${y'_0}^{(d-1)}\cup {z'_0}^{(d-1)}$ and
${P'_0}{\diagdown}{({Y'_0}^{(d-1)}\cup {Z'_0}^{(d-1)})}$
is positive.

Denote by ${x'_0}^{(d-2)}$ the simplex ${y'_0}^{(d-1)}\cap {X'_0}^{(d-2)}$. 
Obviously, ${x'_0}^{(d-2)}\in \xi'$.

\begin{figure}
\begin{center}
\includegraphics[width=0.5\textwidth]{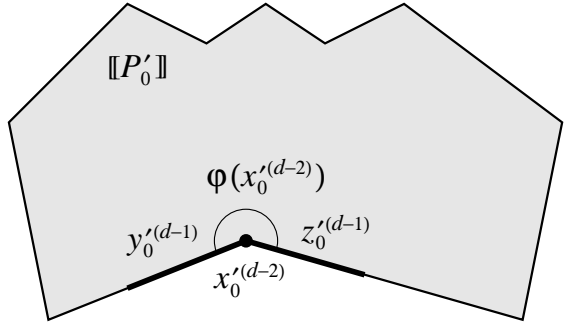}
\end{center}
\caption{The cross-section of the domain
$[\![P'_0]\!]\subset\mathbb{R}^d$ by a two-dimensional plane, 
which passes through an interior point of the simplex ${x'_0}^{(d-2)}$
and which is orthogonal to the affine hull of ${x'_0}^{(d-2)}$.
Intersections with the simplicies ${x'_0}^{(d-2)}$, ${y'_0}^{(d-1)}$ and ${z'_0}^{(d-1)}$ 
are shown in bold}\label{fig1}
\end{figure}

In Figure \ref{fig1} we show schematically the cross-section of the domain
$[\![P'_0]\!]\subset\mathbb{R}^d$ by a two-dimensional plane, 
which passes through an interior point of the simplex ${x'_0}^{(d-2)}$
and which is orthogonal to the affine hull of this simplex.

\textit{Step 4.} 
Choose continuous families $\{v'_s\}_{s\in[0,1]}$ and  $\{w'_s\}_{s\in[0,1]}$ 
of points in $\mathbb{R}^d$ such that the following properties hold true
for all $s\in[0,1]$:

$(\sigma_1)$
$v'_s\in [\![P'_s]\!]$ and $w'_s\in [\![P'_s]\!]$;

$(\sigma_2)$
there exists an isometry $f_s:\mathbb{R}^d\to \mathbb{R}^d$
such that $f_s({y'_0}^{(d-1)})={y'_s}^{(d-1)}$ and $f_s(v'_0)=v'_s$;

$(\sigma_3)$
there exists an isometry $g_s:\mathbb{R}^d\to \mathbb{R}^d$
such that $g_s({z'_0}^{(d-1)})={z'_s}^{(d-1)}$ and
$g_s(w'_0)=w'_s$;

$(\sigma_4)$
the set $\partial(\mbox{conv}\,({y'_s}^{(d-1)}\cup\{v'_s\})){\diagdown}{y'_s}^{(d-1)}$
has no common points with $P'_s$;

$(\sigma_5)$
the set $\partial(\mbox{conv}\,({z'_s}^{(d-1)}\cup\{w'_s\})){\diagdown}{z'_s}^{(d-1)}$
has no common points with $P'_s$;

$(\sigma_6)$
the sum of the internal dihedral angles of the simplicies
$\mbox{conv}\,({y'_s}^{(d-1)}\cup\{v'_s\})$ and $\mbox{conv}\,({z'_s}^{(d-1)}\cup\{w'_s\})$
at ${x'_s}^{(d-2)}$ is equal to $\varphi_*$, where $\varphi_*$ was determined by the formula (\ref{eq6});

$(\sigma_7)$
the Hausdorff distance between the sets $P'_0\cup\mbox{conv}\,({y'_0}^{(d-1)}\cup\{v'_0\})\cup
\mbox{conv}\,({z'_0}^{(d-1)}\cup\{w'_0\})$ and $P'_0$ is less than $\varepsilon/3$.

Obviously, the property $(\sigma_6)$ can be satisfied for $s=0$ because, 
by a suitable choice of the points $v'_0$ and $w'_0$,
the dihedral angle at ${x'_0}^{(d-2)}$ in each of the simplicies
$\mbox{conv}\, ({y'_0}^{(d-1)}\cup\{v'_0\})$ and
$\mbox{conv}\, ({z'_0}^{(d-1)}\cup\{w'_0\})$
can be made equal to any number in the interval $(0,\pi)$, 
while $0<\varphi_*<2\pi$.

It is just as simple to verify that the properties $(\sigma_4)$, $(\sigma_5)$, and $(\sigma_7)$
can also be fulfilled for $s=0$.
Hence, there exists $\lambda>0$ such that the properties $(\sigma_4)-(\sigma_7)$ are fulfilled for 
all $s\in[0,\lambda]$.
Substituting $s\mapsto s/\lambda$ , we may assume without loss of generality that 
$\lambda =1$, i.\,e., that the properties $(\sigma_1)-(\sigma_7)$ are valid for all $s\in [0,1]$.

\begin{figure}
\begin{center}
\includegraphics[width=0.5\textwidth]{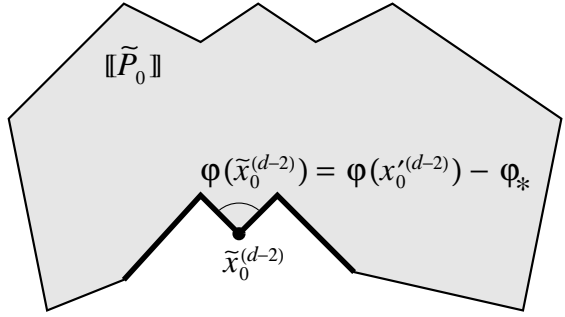}
\end{center}
\caption{The cross-section of the domain
$[\![\widetilde{P}_0]\!]\subset\mathbb{R}^d$ by a two-dimensional plane, 
which passes through an interior point of the simplex ${\widetilde{x}_0}^{(d-2)}$
and which is orthogonal to the the affine hull of ${\widetilde{x}_0}^{(d-2)}$.
Intersections with ${\widetilde{x}_0}^{(d-2)}$ and the sets
$\partial(\mbox{conv}\,({y'_0}^{(d-1)}\cup\{v'_0\})){\diagdown}{y'_0}^{(d-1)}$ and
$\partial(\mbox{conv}\,({z'_0}^{(d-1)}\cup\{w'_0\})){\diagdown}{z'_0}^{(d-1)}$
are shown in bold}\label{fig2}
\end{figure}

\textit{Step 5.}
Let $\mbox{relint}\,({y'_s}^{(d-1)})$ denote the relative interior of the simplex ${y'_s}^{(d-1)}$,
i.\,e., the interior of ${y'_s}^{(d-1)}$ within its affine hull.
For every $s\in[0,1)$, replace the simplicies ${y'_s}^{(d-1)}$ and ${z'_s}^{(d-1)}$ in $P'_s$
by the polyhedra 
$\partial(\mbox{conv}\,({y'_s}^{(d-1)}\cup\{v'_s\})){\diagdown}\mbox{relint}\,({y'_s}^{(d-1)})$ 
and 
$\partial(\mbox{conv}\,({z'_s}^{(d-1)}\cup\{w'_s\})){\diagdown}\mbox{relint}\,({z'_s}^{(d-1)})$
respectively (see Figure \ref{fig2}).
Denote the resulting polyhedron by $\widetilde{P}_s$.

Immediately from the construction of $\widetilde{P}_s$ it follows that the polyhedron $\widetilde{P}_0$ 
and family $\widetilde{\mathscr {F}}=\{\widetilde{P}_s\}_{s\in[0,1)}$ possess the properties
$(\widetilde{\alpha})-(\widetilde{\gamma})$.
\hfill$\square$

\section{Proof of Theorem \ref{thrm1}}\label{s4}

\textbf{Proof} \  Let $P_0$ be an embedded flexible polyhedron in $\mathbb{R}^d$, $d\geqslant 3$, 
and let $\varepsilon>0$. According to Lemma \ref{l2}, there exists an embedded flexible polyhedron 
$\widetilde{P}_0\subset\mathbb{R}^d$ and its flex $\mathscr{\widetilde{F}}=\{\widetilde{P}_s\}_{s\in [0,1)}$
such that the properties $(\widetilde{\alpha})-(\widetilde{\gamma})$ are fulfilled. 

For $j=1,2$ and $s\in[0,1)$, denote by $\mathscr{N}_s^{(j)}(k)$ 
the eigenvalue counting function of the problem (\ref{eq1}) (for $j=1$) or problem (\ref{eq2}) (for $j=2$)
in the domain $\Omega=[\![\widetilde{P}_s]\!]$.

Putting $p=2$ and applying Lemma \ref{l1} for $[\![\widetilde{P}_s]\!]$, we get 
\begin{equation}\label{eq7}
\int\limits_0^k(k-t)^2\,d\mathscr{N}_s^{(j)}(t)=
\frac{2a_d^s}{(d+1)(d+2)}k^{d+2}+
\frac{2a_{d-1}^s}{d(d+1)}k^{d+1}+ 
\frac{2a_{d-2}^s}{(d-1)d}k^d+ 
O(k^{d-1})
\end{equation}
as $k\to\infty$, where
\begin{eqnarray}
a_d^s&=&\frac{\mbox{\rm{vol}}_d([\![\widetilde{P}_s]\!])}{2^d \pi^{d/2} 
\Gamma\bigl(\frac{d}{2}+1\bigr)},\nonumber\\
a_{d-1}^s&=&\mp
\frac{\mbox{\rm{vol}}_{d-1}(\widetilde{P}_s)}{2^{d+1}\pi^{(d-1)/2}
\Gamma\bigl(\frac{d+1}{2}\bigr)},\label{eq8}\\
a_{d-2}^s&=&\frac{1}{2^{d+1}\pi^{d/2}\Gamma\bigl(\frac{d}{2}\bigr)}
\sum\limits_{\widetilde{x}^{(d-2)}_s}\frac{[\varphi(\widetilde{x}^{(d-2)}_s)]^2-\pi^2}
{3\varphi(\widetilde{x}^{(d-2)}_s)}\mbox{\rm{vol}}_{d-2}
\bigl(\widetilde{x}^{(d-2)}_s\bigr).\label{eq9}
\end{eqnarray}
In the formula (\ref{eq8}), the minus sign corresponds to the problem (\ref{eq1}),
while the plus sign corresponds to the problem (\ref{eq2}).
In the formula (\ref{eq9}), the sum is taken over all $(d-2)$-dimensional faces
$\widetilde{x}^{(d-2)}_s$ of $\widetilde{P}_s$ and $\varphi(\widetilde{x}^{(d-2)}_s)$
stands for the value of the dihedral angle of $[\![\widetilde{P}_s]\!]$
at $\widetilde{x}^{(d-2)}_s$.

Assuming that the Dirichlet or Neumann spectrum of the Laplacian in $[\![\widetilde{P}_s]\!]$
does not change when $s$ changes in the interval $[0,1)$, we get that, for $j=1$ or $j=2$, 
the eigenvalue counting function $\mathscr{N}_s^{(j)}(k)$ is independent of $s$.
Hence, for $j=1$ or $j=2$, the left-hand side of the formula (\ref{eq7}) is independent of $s$, i.\,e., 
the formula (\ref{eq7}) is an asymptotic expansion of the function
\begin{equation*}
k\mapsto\int\limits_0^k(k-t)^2\,d\mathscr{N}_s^{(j)}(t)
\end{equation*}
which is independent of $s$. Since the asymptotic expansion of a function is determined uniquely,
it follows that the coefficients at $k^{d+2}$, $k^{d+1}$ and $k^d$ on the right-hand side of the 
formula (\ref{eq7}) are also independent of $s$.

However, the latter statement is false for the coefficient at $k^d$, since $a_{d-2}^s$ 
is not constant with respect to $s$. To prove this, we recall that, according to Lemma \ref{l2},
there exists a $(d-2)$-dimensional face $\widetilde{x}_0$ of $\widetilde{P}_0$ such that,
for the corresponding $(d-2)$-dimensional face $\widetilde{x}_s$ of $\widetilde{P}_s$, we have

$\bullet$ $\varphi(\widetilde{x} _s)\to 0$ as $s\to 1$; and

$\bullet$ the value of the dihedral angle of $[\![\widetilde{P}_s]\!]$ at every
$(d-2)$-dimensional face $\widetilde{y}_s$ of $\widetilde{P}_s$, different from 
$\widetilde{x}_s$, is uniformly separated from 0 and $2\pi$: 
$\widetilde{m}<\varphi(\widetilde{y}_s)<\widetilde{M}$.

Therefore, in the sum on the right-hand side of the formula (\ref{eq9}), exactly one term 
(namely, the term corresponding to the face $\widetilde{x}^{(d-2)}_s=\widetilde{x}_s$),
tends to infinity as $s\to 1$, while each of the remaining terms is uniformly bounded for $s\in [0,1)$.

Hence, choosing $s$ sufficiently close to 1, we can to make $a_{d-2}^s$ as large as we want. 
Thus, $a_{d-2}^s$ is not constant in $s$, and both Dirichlet and Neumann spectra of the 
$d$-dimensional Laplace operator in the domain $[\![\widetilde{P}_s]\!]\subset\mathbb{R}^d$ 
do not remain unaltered when $s$ changes in the interval $[0,1)$.
\hfill $\square$

\textbf{Remark 3} \ 
Observe that the coefficients $a_d^s$ and $a_{d-1}^s$ in the formula (\ref{eq7}) are constant in $s$. 
This fact underlines the non-triviality of the conjecture on the invariance of the Dirichlet and 
Neumann spectra of the Laplacian  during the flex of the boundary of a domain, flipped by Theorem \ref{thrm1}.

\section*{Acknowledgements}\label{s5}
The author is grateful to Dr. Evgeni\u{\i} P. Volokitin for assistance in preparation of the figures and
to Prof. Alexey Yu. Kokotov for bringing his attention to \cite{Wa05}.

\noindent{Victor Alexandrov}

\noindent\textit{Sobolev Institute of Mathematics}

\noindent\textit{Koptyug ave., 4}

\noindent\textit{Novosibirsk, 630090, Russia}

and

\noindent\textit{Department of Physics}

\noindent\textit{Novosibirsk State University}

\noindent\textit{Pirogov str., 2}

\noindent\textit{Novosibirsk, 630090, Russia}

\noindent\textit{e-mail: alex@math.nsc.ru}

\end{document}